\documentclass[preprintnumbers,amsmath,amssymb]{revtex4}
\usepackage{epsfig}
\usepackage{graphicx}
\usepackage{dcolumn}
\usepackage{bm}

\usepackage{amssymb}

\begin{document}



\title{Emergence of Synchronous Oscillations in Neural Networks Excited by Noise.}


\author{M.-P. Zorzano }
\email{zorzanomm@inta.es}
\homepage{http://www.cab.inta.es}
\affiliation{Centro de Astrobiolog\'{\i}a (CSIC-INTA),
        Carretera de Ajalvir km 4,    28850 Torrej\'{o}n de Ardoz, Madrid, Spain}

\author{L. V\'azquez}
\email{lvazquez@fdi.ucm.es}

\affiliation{Departamento de Matem\'{a}tica Aplicada. Facultad de
    Inform\'{a}tica. Universidad Complutense. 28040 Madrid}
\affiliation{Centro de Astrobiolog\'{\i}a (CSIC-INTA),
        Carretera de Ajalvir km 4,    28850 Torrej\'{o}n de Ardoz, Madrid,
    Spain}
\date{\today}

\begin{abstract}

The presence of noise in non linear dynamical systems can play a constructive role,
increasing the degree of order and coherence or evoking improvements in the
performance of the system. An example of this positive influence in a
biological system is the impulse transmission in neurons and the
synchronization of a neural network. Integrating numerically the Fokker-Planck
equation we show a self-induced synchronized oscillation. Such an oscillatory state appears in a
neural network coupled with a feedback term, when this system is
excited by noise and the noise strength is within a certain range.

\end{abstract}

\maketitle


\section{Introduction.}

Neural networks consist of many
nonlinear components (neurons that beyond a threshold emit action potentials) which are interdependent (the output of a neuron is the input
 of another neuron in the network) and form a complex system with new emergent
 properties that are not hold by each individual item in the system alone.
The emergent property of this dynamical system is that a set of
 neurons will synchronize and fire impulses simultaneously. In the context of neuroscience,
 this emergent property is used to implement quite sophisticated and highly specialized "logical"
functionalities  such as memorization with Hebbian learning, and
 recognition of patterns (or memories).

There at least two mechanisms by which such synchronous oscillations can take place:
\begin{itemize}
\item{Synchrony can be a consequence of a common input produced by an
    oscillating neuron (or set of neurons) (pacemakers).}
\item{Synchrony can also be a consequence of an emergent population
    oscillation within a network of cells. There is no external coordination
    when this oscillation is built up, it is self-organized and its
properties must be related to the structural way in which the networks are
connected. } 
\end{itemize}

Each neuron can receive input from  neighboring neurons of the network and from external
 impulses. These neural networks are also excited by electrical noise 
 which is ubiquitous in the neuronal system and seem to be able to operate in
 such a noisy environment in a
 robust way.  The
 main sources of noise are related to the synaptic connections and
 voltage-gated channels
 \cite{NoiseIChannel,NoiseSynaptic,AnaCalcium,ChanNoise}. The role of noise in the
 functioning of the non-linear nervous system is poorly understood but there are evidences of  positive interactions of noise and nonlinearity in
 neuronal systems \cite{HumanSR1,paddle,HumanSR2,HumanSR3,RatSR,RatColor}.

In the
present work we will investigate the relevance of  noise in the
synchronization of a neural network.
In order to do so we will integrate numerically the Fokker-Planck equation associated to the
stochastic system. This has the advantage of giving better accuracy in the
tails of the distribution than solving the stochastic differential equations.  Integrating numerically this equation
one can obtain a quick overview of the system dynamics and the time evolution of
the probability density. We will show that noise can help the system building
a synchronous oscillation.

The starting point is the equation of motion of the single
neuron, see section \ref{1}. Then, in section \ref{2}, we describe the ensemble of neurons excited by noise
 using the Fokker-Planck equation, which we integrate numerically in different
 scenarios. We first
 show a noise induced transition in the distribution for a certain range of
 noise intensity (see subsection \ref{2-1}). Next we introduce different situations where synchrony is
attained in the presence of noise: first exciting the ensemble of neurons
with an oscillating external input (see subsection \ref{2-2}), then
introducing only a feedback term (see subsection \ref{2-3}) and finally reinforcing the
external excitation with the internal feedback (see subsection \ref{2-4}). In all these scenarios there exist a
defined range of noise intensity where the
synchronization is maximized.

\section{FitzHugh-Nagumo equation for a single neuron.}\label{1}

The impulse transmission in a neuron can be essentially
described with the Hodgkin-Huxley model, \cite{HH}. The
FitzHugh-Nagumo (FHN) model (\cite{Fitz,Nagumo}) is a simplified variant of the previous model
accounting for the essentials of the regenerative firing mechanism in an
excitable nerve cell, namely, it has a stable rest state,
and with an adequate amount of disturbance it generates a pulse with a
characteristic magnitude of height and width.

 In the FHN model all dynamical variables of the
neuron are reduced to two quantities: the voltage of the membrane $u$ (fast
variable) and the recovery variable $v$ which corresponds to the refractory properties of the membrane (slow variable):
\begin{eqnarray}
\dot{u}&=&c(-v+u-u^{3}/3 + I_1(t)+\sqrt{2D}\xi (t))\\
\dot{v}&=&u-bv+a\\
<\xi(t)\xi(t')>&=&\delta(t-t')
\end{eqnarray}
The neuron is excited by an external input $I_t=I_1(t)+\sqrt{2D}\xi (t)$. Here we have added white
Gaussian noise
$\xi(t)$ to the
standard FHN equation obtaining a Stochastic Differential Equation (SDE).  
For all the cases described below we will study this system for
$c=10,a=0.7$ and $b=0.8$.  The neuron is said to "fire" or emit a "spike" when
$u>0$. This
 impulse is then transmitted along the axon and elicits the
 emission of neurotransmitters which in turn excite other neurons.

\begin{figure}[h]
\begin{center}
\epsfig{figure=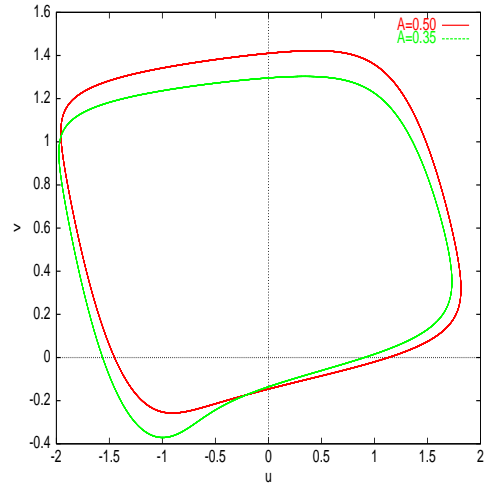,height=6.5cm,width=6.5cm}
\caption{\label{detFHN}Phase space projection of two supra-threshold
  trajectories for the FHN excited with constant amplitude $I_t=A$. 
}

\end{center}
\end{figure}

Let us first consider the deterministic system ($D=0$) with constant excitation
 $I_t=A$. A typical feature of a neuron system 
 is that the information is encoded in the firing rate. For constant input there is a threshold amplitude (in this case
 $A=0.35$) beyond which the output is a periodic
 oscillating system (a train of spikes with $u>0$), whose frequency (and not
 the amplitude)
 depends on the amplitude $A$ of the exciting force.   In Fig. \ref{detFHN}-left we show the
 projection onto the phase space of two of these attractors excited by
 different inputs of constant amplitude. All the
 orbits are attracted to this periodic attractor.  Such attractors will also
 be seen later in the extension to a noisy system. If the excitation is below
 the threshold then the system goes to the stable rest state (there is no
 oscillation in this case). 

\begin{figure}[h]
\begin{center}
\epsfig{figure=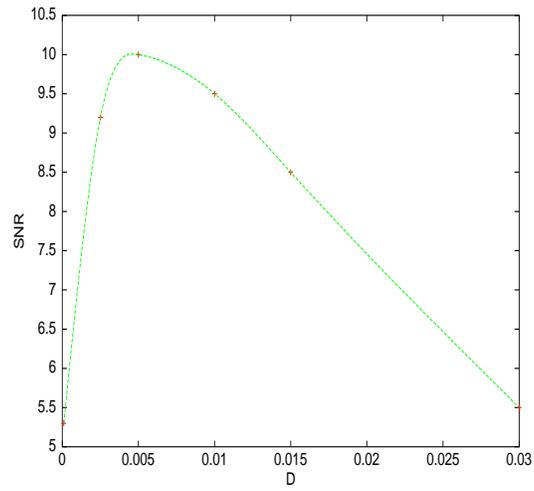,height=6.5cm,width=7cm}
\caption{\label{stoFHN}The signal to noise ratio shows the presence of
  stochastic resonance in the SDE of a FHN system excited by
  $I_t=A\cos(2\pi f t)+\sqrt{2D}\xi(t)$ when $A=0.17$ and $f=0.55$.}
\end{center}
\end{figure}

We study now the response of the
system to a time varying force of the type $I_t=A \cos{(2\pi f
  t)}+\sqrt{2D}\xi(t)$, with $A=0.15$, $f=0.55$ and various noise
intensities $D$. For a given initial condition ($u(t_0),v(t_0)$)=(0,0), we perform a Monte Carlo
simulation: integrate the SDE for
different noise realizations
and average over the ensemble to obtain $<u(t)>$. We can evaluate the
spectrum and the Signal-to-Noise-Ratio (SNR). The system shows a maximal
susceptibility to the external periodic forcing for $D\approx 0.005$, see
Fig. \ref{stoFHN}.  This feature will also survive
when we consider an ensemble of neurons and will be decisive for the behavior
of the network.

 The enhancement of the system response due to noise is a characteristic trace of stochastic resonance  (for a comprehensive
 review on prominent references see \cite{GammaHaengi,NeimanSR}). Noise can induce hopping from the
resting state to excited state and vice-versa. The mean time that it takes a
neuron to move from one state to the other for a given noise intensity is
 called mean scape time. Stochastic resonance takes
place when there exists some kind of
synchronization between the forcing time and the mean scape time, i.e. when
 there exists a time-scale resonance condition.  For
 instance, for the double-well
potential, this happens when the
mean scape time is comparable with  half
the period of the periodic
forcing.


\clearpage
\newpage
\section{Neuron ensemble and Fokker-Planck equation.}\label{2}

In the previous section we have introduced the behavior of a single neuron
described by an FHN oscillator excited by constant force or by a periodic
force plus noise. But what happens when we consider an ensemble of such
excitable elements?.
In order to describe this set of oscillators we will use the
Fokker-Planck (FP) equation, which is a partial differential equation, describing the
time evolution of the probability density associated to a SDE. If one needs to obtain global information from a system of
independent neurons excited by noise, one can either
run a Monte Carlo simulation for the SDE for all the initial conditions (one
per neuron of the network) and group the solution in histograms or
solve the FP equation for the density in phase space with a given initial
condition \cite{MpMaisLV}. Integrating the FP equation we get better accuracy in
the tails, where the probability is very low, than running the Monte Carlo
simulation. These tails will be particularly important in the FHN
stochastic equation since the supra-threshold excursions in the attractor
region, that are relevant for the network excitation, are of low probability.

We integrate numerically the Fokker-Planck equation of the previously described FHN
system: 
\begin{eqnarray}
\frac{\partial \rho}{\partial t}&=&-\frac{\partial(u-bv+a) \rho}{\partial
  v}-\frac{\partial c(-v+u-u^{3}/3 + I_1(t)) \rho}{\partial u}+Dc^{2}\frac{\partial^{2} \rho}{\partial u^{2}}\\
\end{eqnarray}
using
an alternating semi-implicit scheme, to obtain the time evolution of the probability density $\rho(x,v,t_0)$ for any
given initial condition. This equation gives the
distribution of the state of our network of neurons, the neurons can
 be excited or at rest depending on the input current and
noise. If we monitor the activity of this network from outside, as it is done
when performing an electro-encephalograph, we measure the average
response of this network $<u>$. We will see that this value can be at rest or oscillate with
different frequencies depending on the input current
$I(t)=I_1(t)+\sqrt{2D}\xi(t)$, with $\xi(t)$ white Gaussian noise. The
connectivity of the network is implicitly given in the structure of $I_1(t)$.

For the rest of the work we will integrate this partial differential equation with the following values: $\Delta
u=0.03$, $\Delta v=0.013$, absorbing boundaries at $u=\pm 4,5,v=\pm 2.34$ and
$\Delta t=0.01$  for various noise
intensities $D$. Our initial condition is a Gaussian
 distribution centered close to the equilibrium  point (non firing condition),
 representing a network with all the neurons at rest. The Gaussian initial
 condition has mean values $<u>=-1.0$
and $<v>=-0.55$ and variances $\sigma_u^{2}=0.05, \sigma_v^{2}=0.013$.

\subsection{Neural network excited only by noise.}\label{2-1}

\begin{figure}[h]
\begin{center}
\epsfig{figure=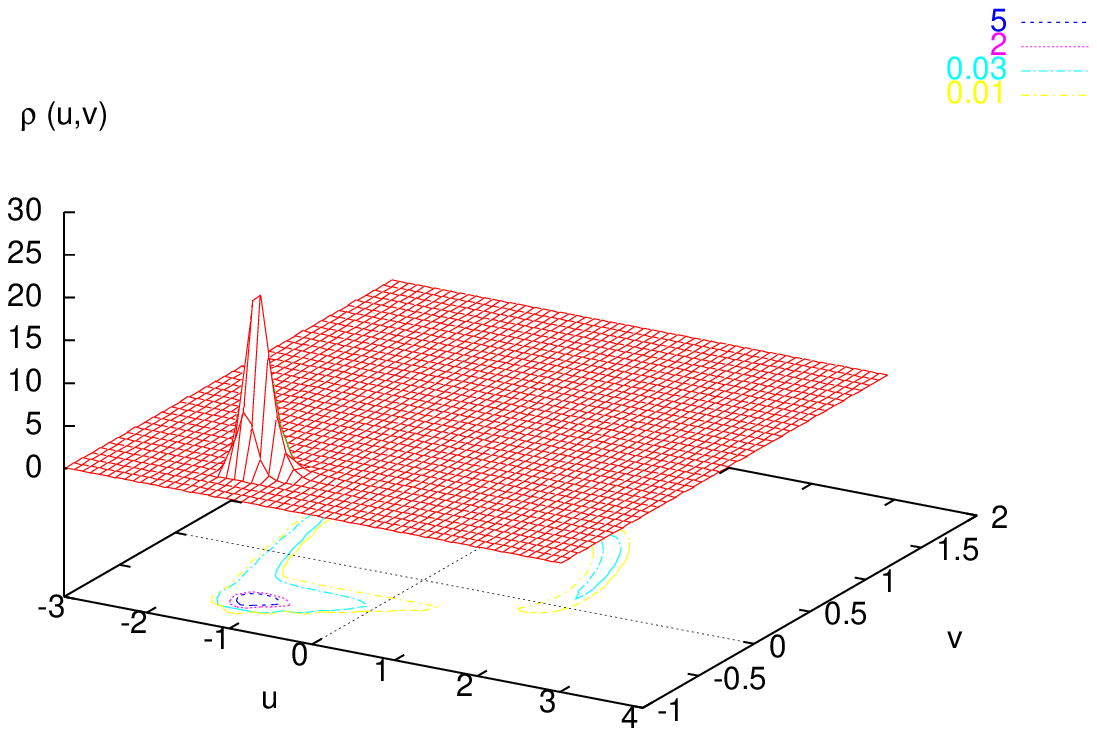,height=7cm,width=6.5cm}
\epsfig{figure=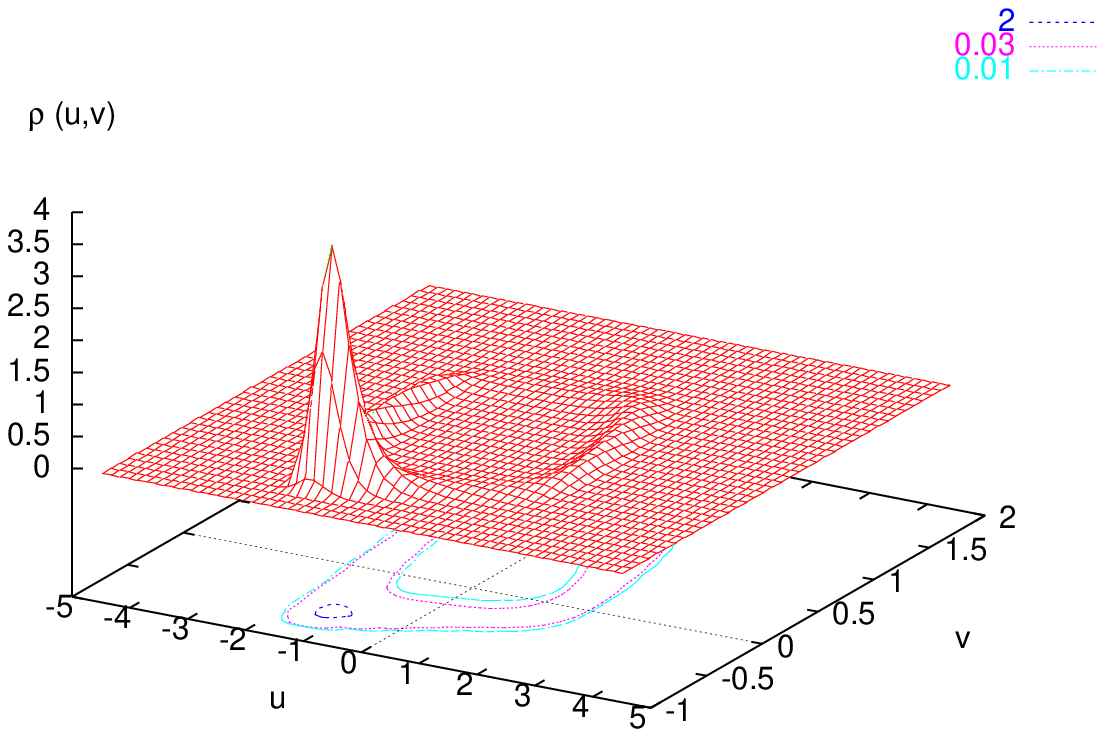,height=7cm,width=6.5cm}
\caption{\label{FHN_FP_A=0D_005}Stationary distribution of the neural network
  distribution excited only by noise with intensity
  $D=0.001$ (left) and $D=0.005$ (right). The tails of the distribution extend
  to the supra-threshold area,
  lying on top of the attractor. There is a significant increment in the
  probability 
  for a neuron to be excited (supra-threshold area) when $D=0.005$.}
\end{center}
\end{figure}

Let us first consider the ensemble of independent neurons when the only input current is the
 noisy term $I_2(t)=\sqrt{2D}\xi(t)$. For low noise intensities,$D=0.001$, after a transient time, the probability of these large
excursions is almost non-existent. See Fig. \ref{FHN_FP_A=0D_005}-left where the density is
peaked on top of the stable fixed point.
For bigger noise intensities, $D=0.005$, after the transient time, there is a certain probability for a neuron to be
excited, even in the absence of external forcing. In Fig. \ref{FHN_FP_A=0D_005}-right we
show the stationary
distribution for $D=0.005$ that has non-vanishing tails on the supra-threshold area.  We can talk
of noise induced transition  in the density distribution in the sense that the stationary
distribution changes with noise. We will see that this strong dependence on the noise intensity will
survive if we take into account additional influences, such as external driving
or feedback, and will be determinant for tuning the system response.

\clearpage
\newpage

\subsection{Periodic force excitation.}\label{2-2}

\begin{figure}[h]
\begin{center}
\epsfig{figure=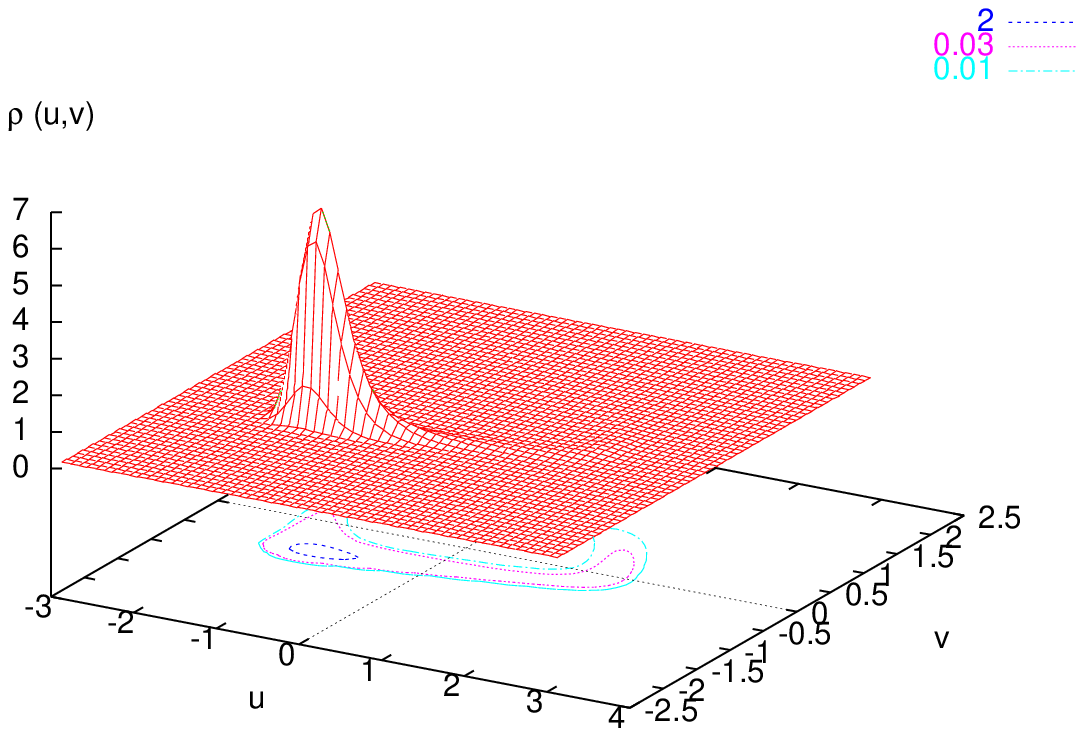,height=7cm,width=6.5cm}
\epsfig{figure=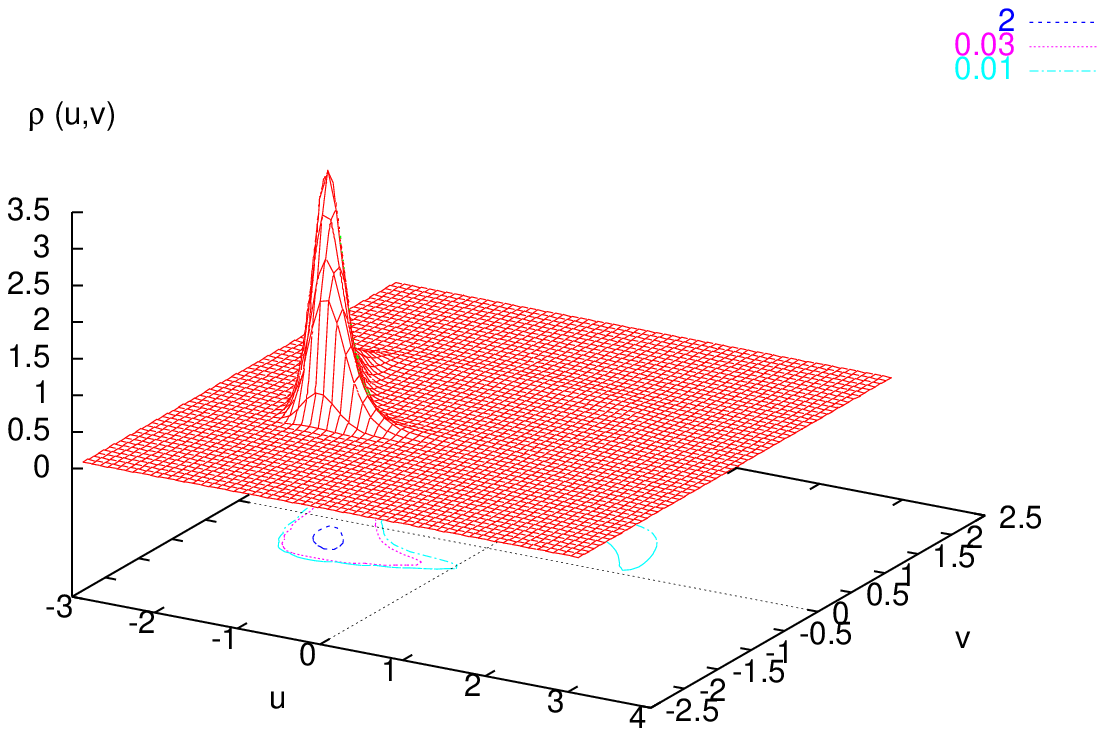,height=7cm,width=6.5cm}
\epsfig{figure=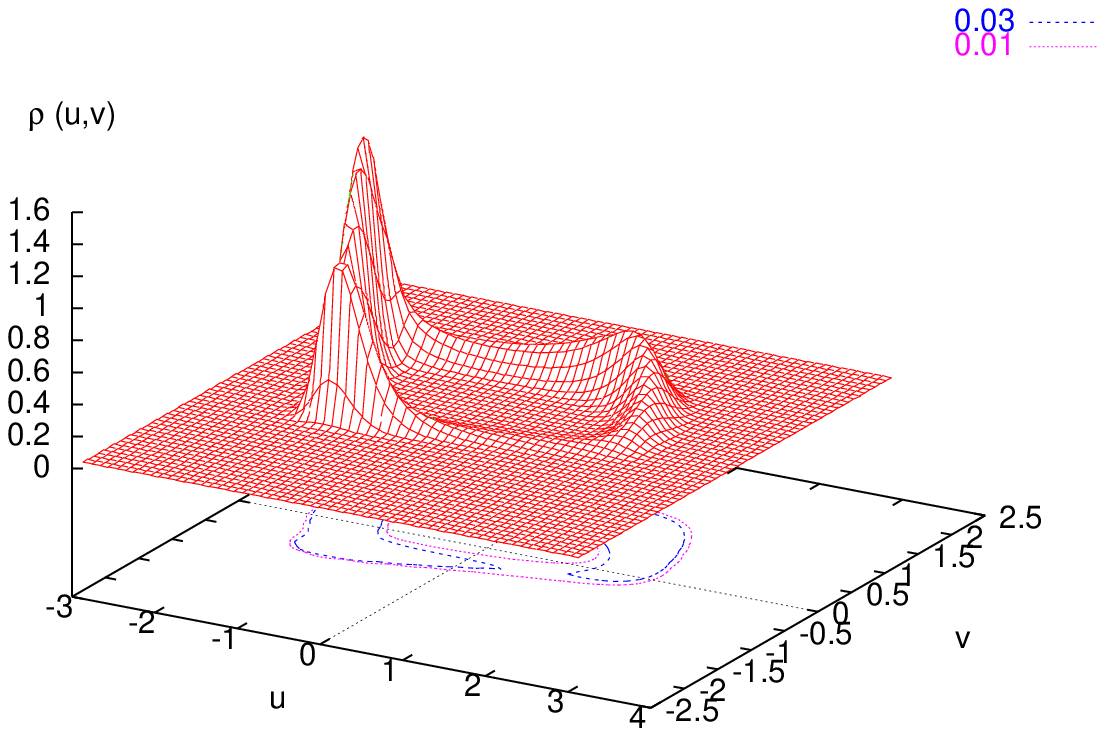,height=7cm,width=6.5cm}
\caption{\label{FHN_FP_A=0.15}Pseudo-stationary distribution (periodic in time) for the neural
  network excited by an external periodic force $I_1(t)=A\cos(2\pi 0.55 t)$ and noise. Shown at
  time $t=100,120,140$ for $D=0.005$ and $A=0.15$.}
\end{center}
\end{figure}

\begin{figure}[h]
\begin{center}
\epsfig{figure=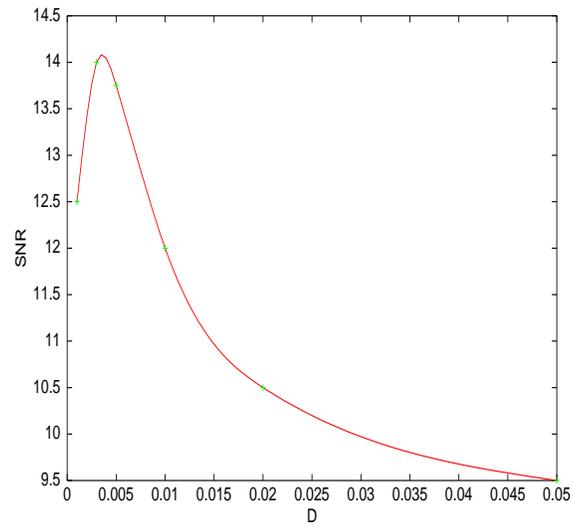,height=7cm,width=7.5cm}
\caption{\label{SNR_A_15_D}Signal to noise ratio from the
  spectral density for $I_t=A\cos(2\pi 0.55 t)$ with $A=0.15$ and various
  noise intensities. The response of the system to the forcing is maximized at
  $D=0.005$.}
\end{center}
\end{figure}

Next we study the system excited by noise for $D=0.005$, and an external
periodic force $I_1(t)=A\cos(2\pi f t)$ (induced for instance by an external pacemaker) with
frequency $f=0.55$ (period $T=1/f$), and $A=0.15$. After the
transient time the density moves along the phase space plane (over a slightly different 
attractor) in a periodic way, with the period of the exciting force, see
Fig. \ref{FHN_FP_A=0.15}. This periodic oscillation is
sustained as long as it is excited by the external periodic force. We have
studied the time evolution of the 1st moment  $<u(t)>=\int_{-\infty}^{\infty}u \rho (u,v,t) du
 dv$ for different noise intensities when the pseudo-stationary solution is
 reached 
As in the single neuron case introduced in section \ref{2-1}, we also find here the stochastic resonance phenomenon in the signal to noise ratio of
the spectral analysis $<u(t)>$, see
 Fig. \ref{SNR_A_15_D}.


\clearpage
\newpage
\subsection{Neural network with excitatory feedback and delay.}\label{2-3}

What happens when our neural network has no external driving force but instead
is coupled internally via some kind of feedback term? We will see that this network can
exhibit autonomous stochastic resonance: for
a certain range of noise intensities a
self-organized and robust coherent oscillation appears in the
absence of external periodic forcing. 

 Let us first find the simplest possible model which includes the basic features of a
 neural network in the brain. 
The current flowing into each neuron is due to the
 interactions with other cells and the response of the neuron is
 expected to depend on the sum of the active synapses (sum of firing neurons where the
 dendrites are connected to). We simulate now this input current $I(t)$ as the
 sum of two terms. The first one $I_1(t)=A\cdot n(t)$ is proportional to the
 number of excited neurons in the network $n(t)=\int_0^{\infty}\rho (u,v,t-\Delta T) du
 dv$, and is the same for each neuron (mean field
 approximation). This term is evaluated at time $t-\Delta T$ but it acts at time
 $t$. We are tacitly assuming that our neural network is fully
 connected, and that there is some transmission delay (synaptic communications
 between neurons depend on propagation of action potentials 
 often over appreciable distances), where the delay $\Delta T$
 corresponds to the time that this process takes before acting back on the
 network. The second term $I_2(t)=\sqrt{2D}\xi(t)$ is stochastic, as above, and
 represents the deviations with respect to this global term for each
 neuron input.  Finally we have a network of
 neurons described as an ensemble of
 FHN oscillators coupled via a nonlocal feedback term.

\begin{figure}[h]
\begin{center}
\epsfig{figure=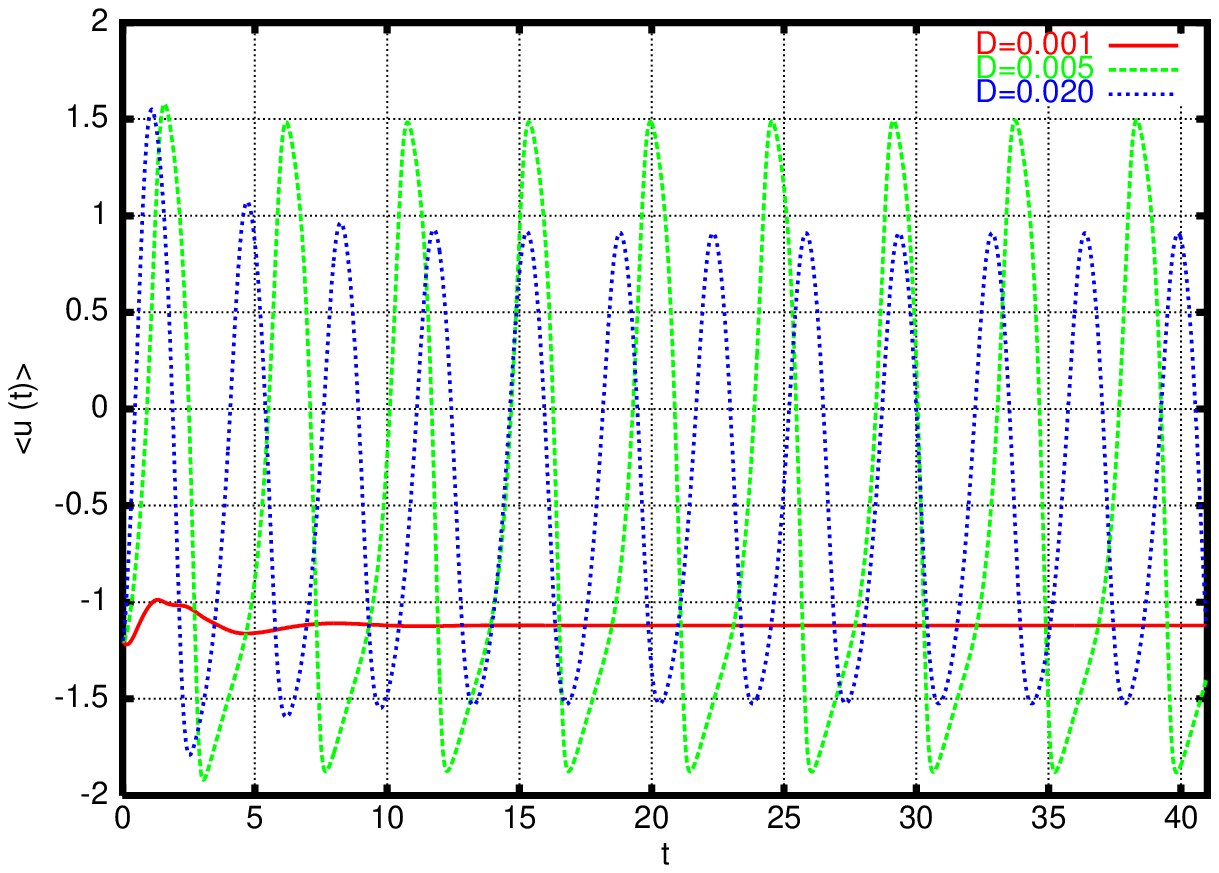,height=6cm,width=6.5cm}
\epsfig{figure=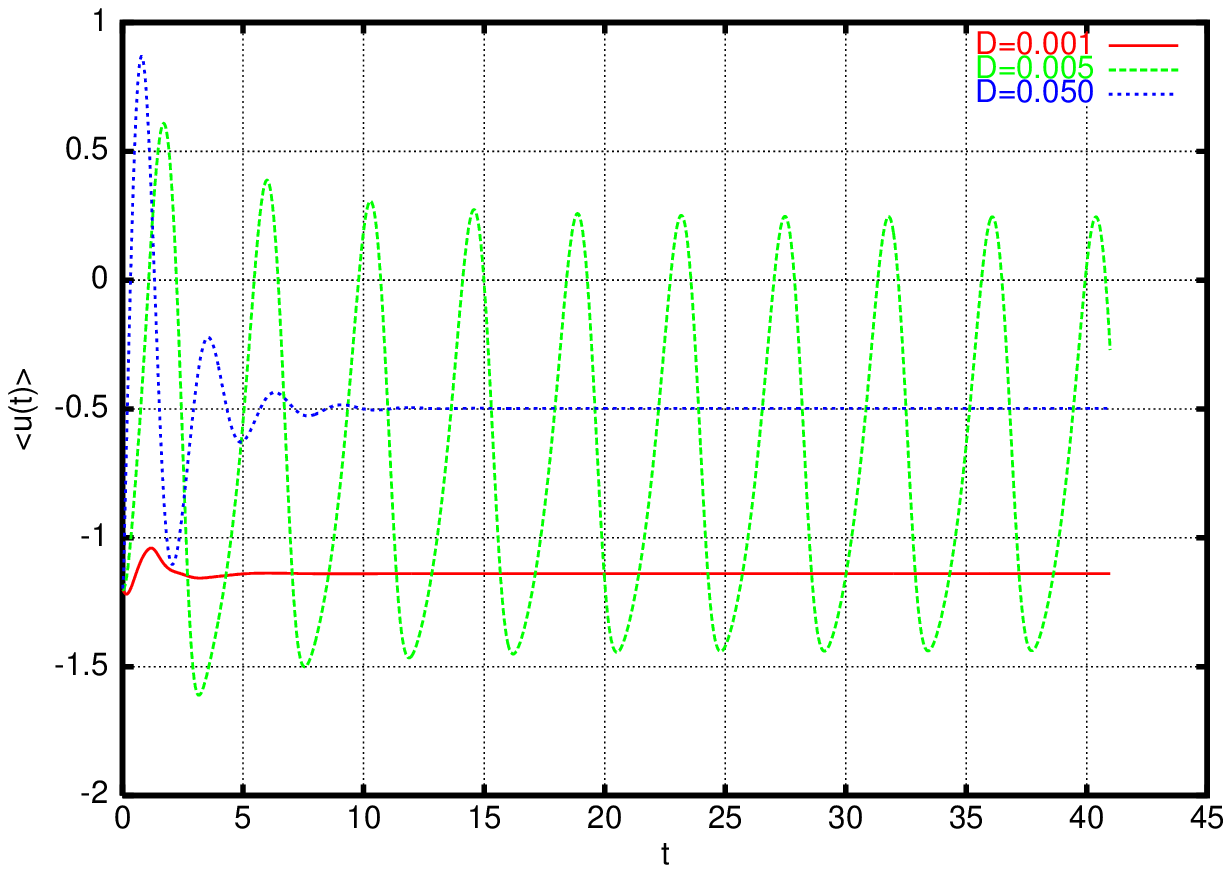,height=6cm,width=6.5cm}
\caption{\label{feedbackFHN}Time evolution of the network output $<u(t)>$ for
  different noise intensities showing and enhancement of the response for
  $D=0.005$, when $A=0.9$ (left) and when $A=0.5$ (right). A self-sustained global
  oscillation is excited by the internal feedback and noise. If the feedback level is
  weak only for noise levels in the region of stochastic resonance ($D=0.005$) can this oscillation emerge.}

\end{center}
\end{figure}

We have mentioned before that the neuron dynamics is such that, for constant input above the threshold of
excitation the system oscillates whereas for input below the threshold the
system is at rest. We first study this system for $\Delta T=0.2$ and a strong feedback term
$A=0.9$. In Fig. \ref{feedbackFHN}-left we show the response of the system $<u(t)>=\int_{-\infty}^{\infty}u
\rho (u,v,t) du dv$ for different noise intensities $D$. For low noise
intensities ($D=0.001$) the feedback term
$I_1(t)$ alone is below the threshold of excitation. The response of the system is enhanced by  stochastic
resonance, reaching the maximal
amplitude and level of synchronization for $D=0.005$. Beyond this noise
intensity the response to the feedback is lower but still beyond the
level of excitation so that a smaller amplitude oscillation is sustained for $D=0.02$. The
maximum value of $n(t)$, $n_{max}$, gives a measurement of the synchronization
of the network, for the cases shown
$n_{max}=0.05$ for $D=0.001$, $n_{max}=0.95$ for $D=0.005$ and $n_{max}=0.8$ for $D=0.02$. 

 If the amplitude of feedback excitation $A$ is weaker, such as
$A=0.5$, the enhancement effect is even clearer, see Fig. \ref{feedbackFHN}-right. The
feedback is only able to overcome the threshold of oscillation in a narrow
region of the noise intensity $D$, where stochastic resonance is at its maximum.

\begin{figure}[h]
\begin{center}
\epsfig{figure=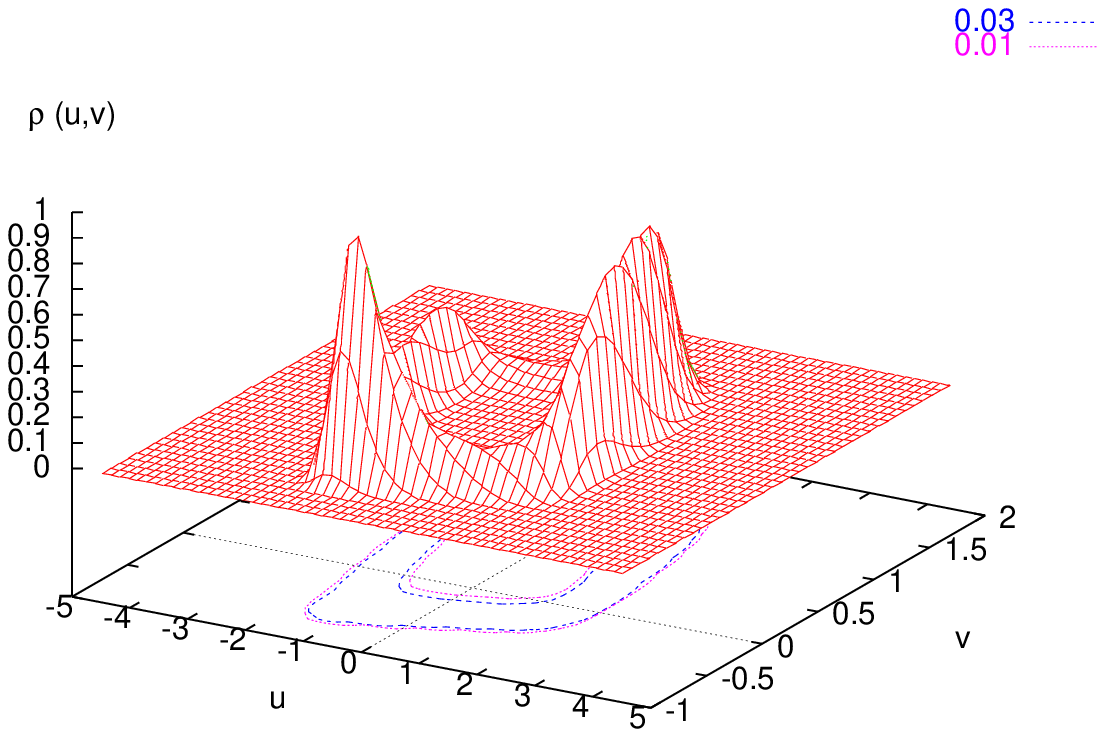,height=7cm,width=6.5cm}
\epsfig{figure=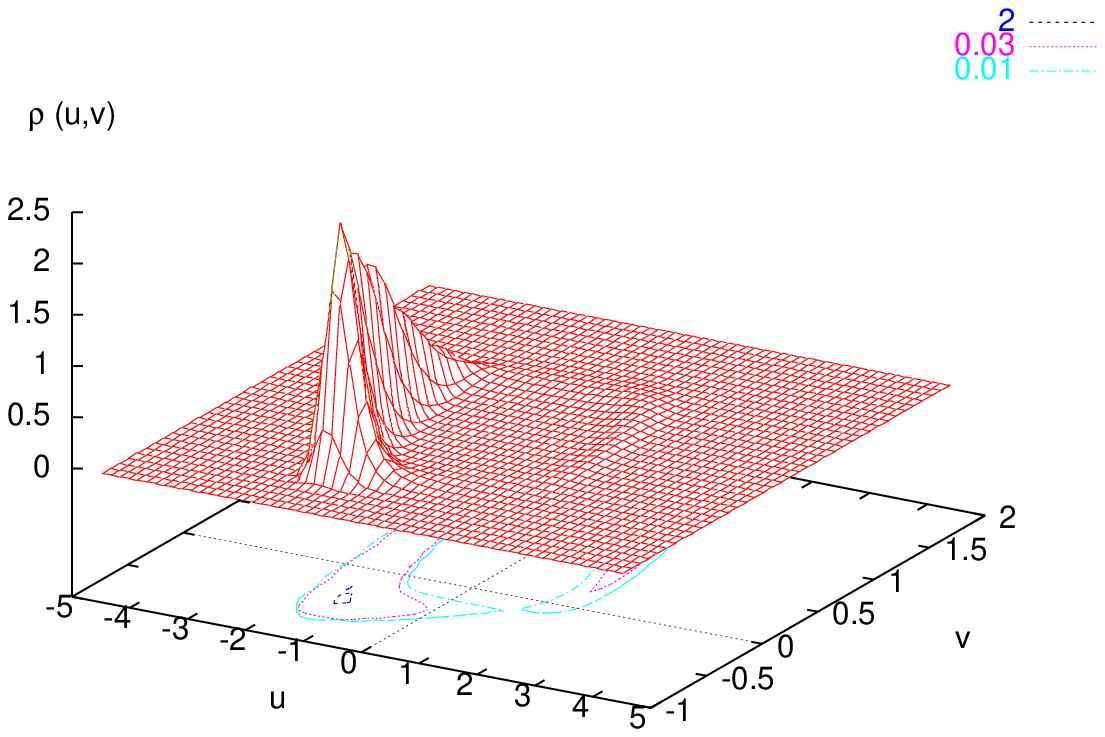,height=7cm,width=6.5cm}
 \epsfig{figure=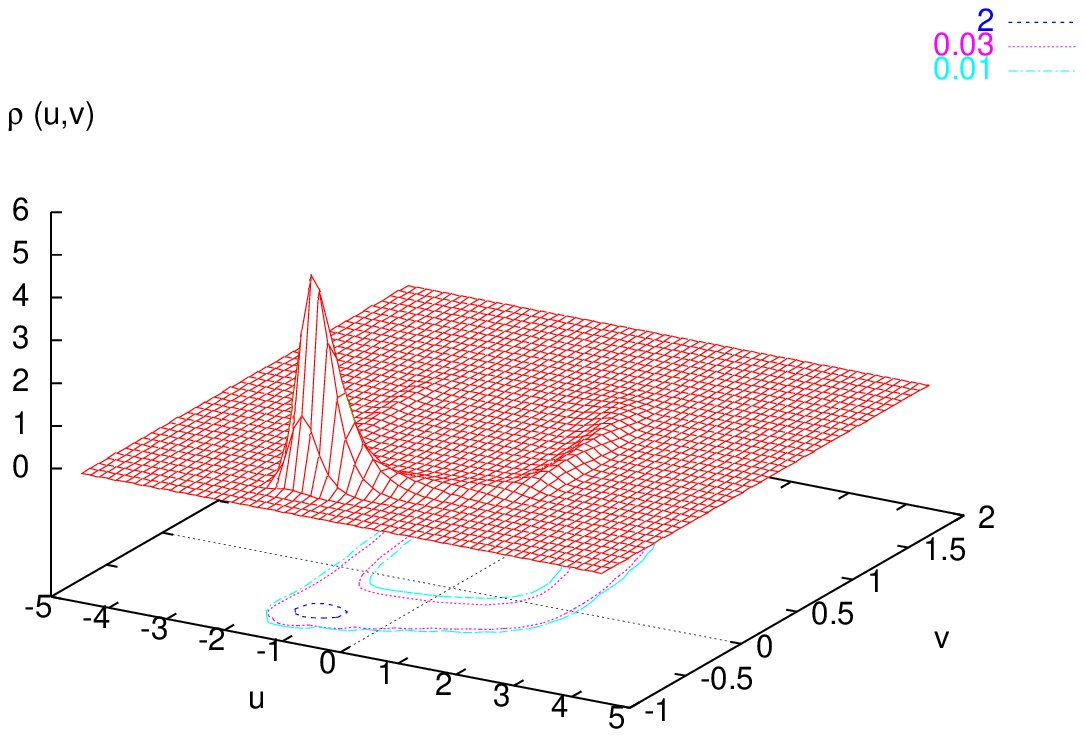,height=7cm,width=6.5cm}
\epsfig{figure=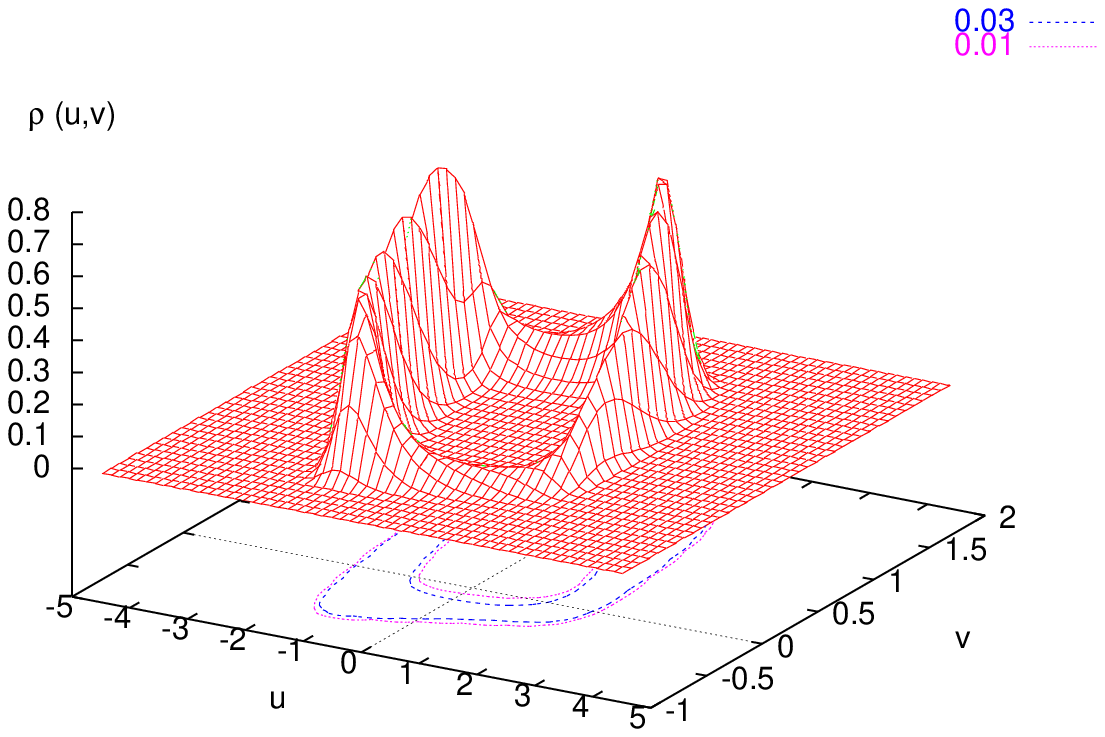,height=7cm,width=6.5cm}
\caption{\label{feedbackFHNrho}Neural network distribution at different times
  when a self-oscillatory state is excited by a weak feedback term ($A=0.5$) and noise ($D=0.005$): t=10.24,20.48,30.72 and 40.96.}
\end{center}
\end{figure}

Notice that  not every single neuron is necessary oscillating
  with the same phase since the density distribution is diffused along the
  attractor, see the density plots at different times in
  Fig. \ref{feedbackFHNrho}.

\begin{figure}[h]
\begin{center}
\epsfig{figure=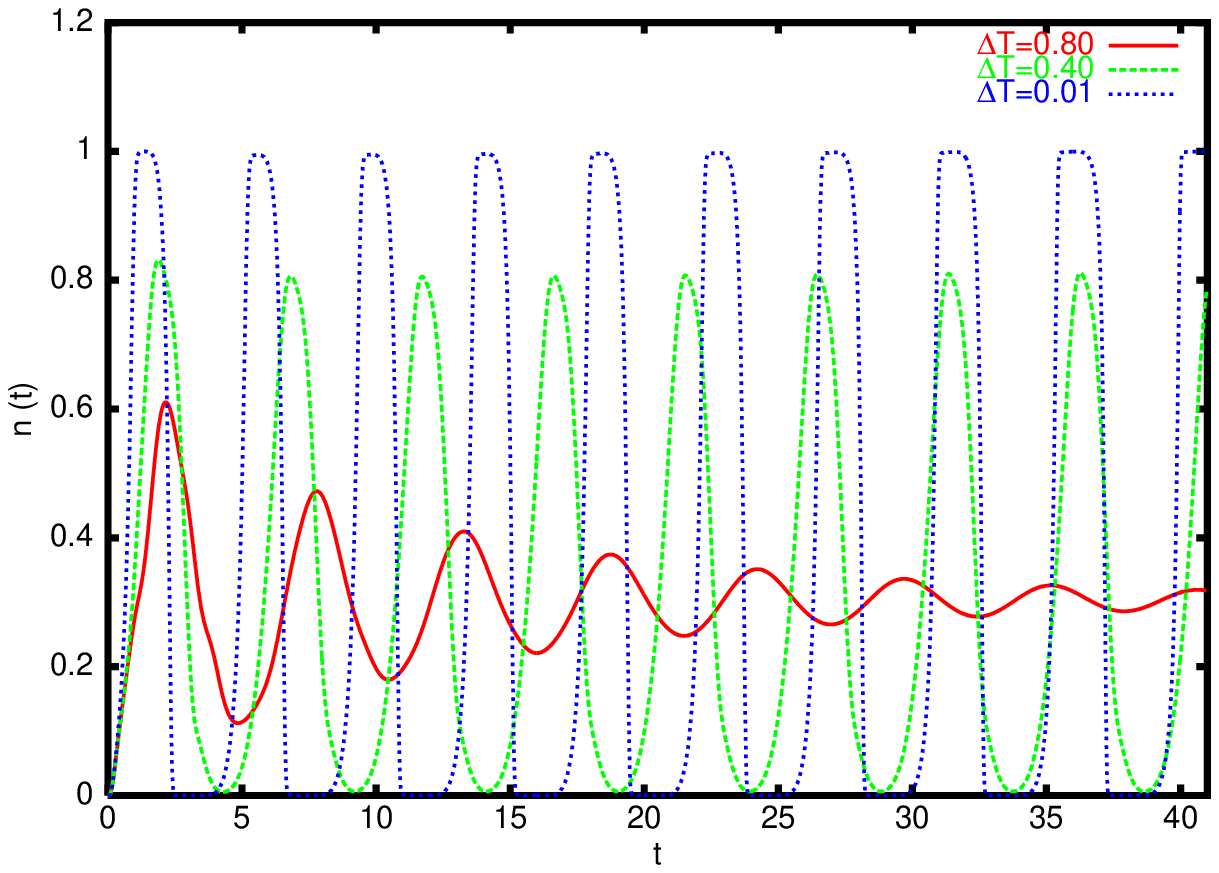,height=6cm,width=6.5cm}
\epsfig{figure=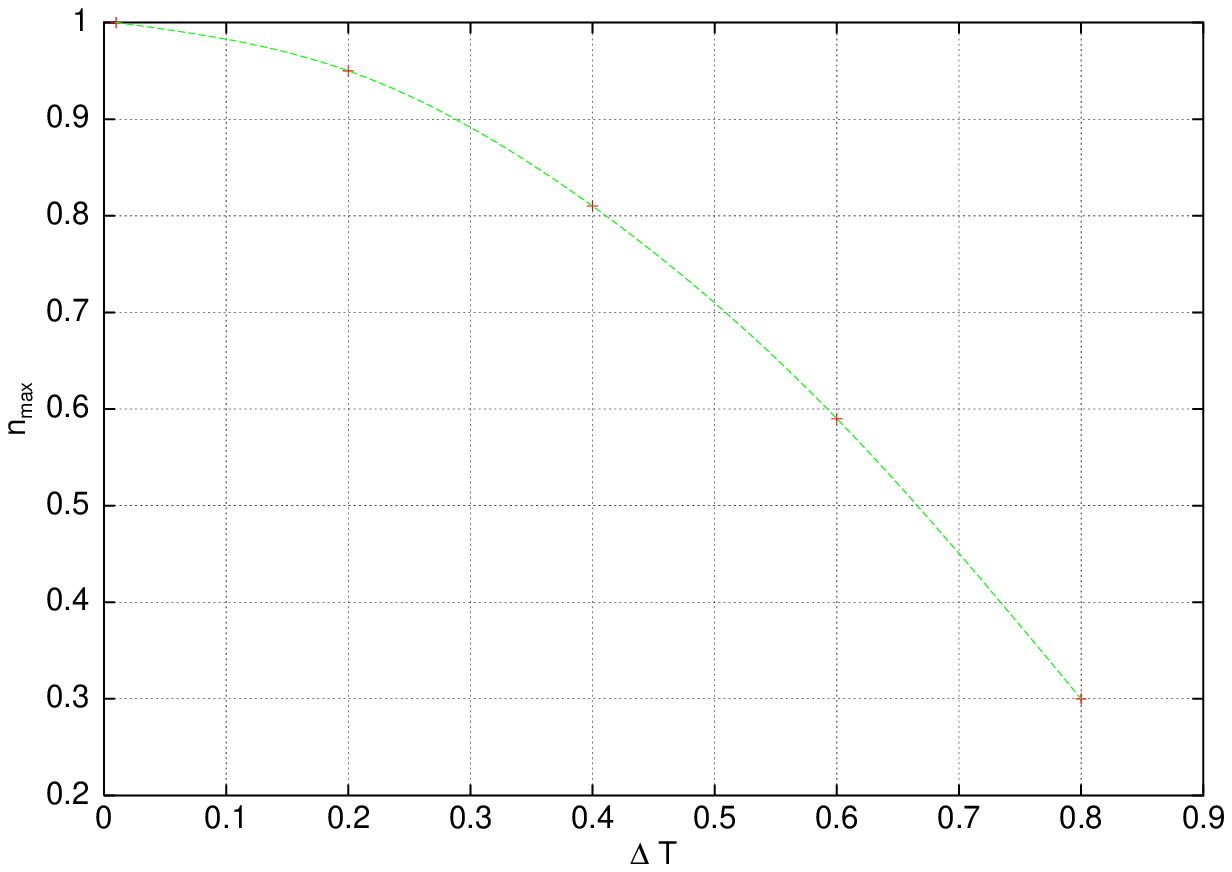,height=6cm,width=6.5cm}
\caption{\label{U_excited_DeltaT} (Left) Fraction of neurons above the threshold $n(t)=\int_0^{\infty}\rho (u,v,t-\Delta T) du
 dv$ as a function of time for different delays $\Delta T$, when the response
 of the system is at its maximum due to stochastic resonance ($D=0.005$). For longer conduction delays the
 synchronization level and the oscillation frequency decrease. (Right) Maximum
 value of $n_{max}$ (synchronization level) as a function of delay $\Delta T$. .}
\end{center}
\end{figure}

Next we keep the noisy input constant and increase the delay time $\Delta
T$ up to
$20\%$ of the neuron cycle. We calculate $n(t)$ the fraction of the population
that is excited at a given time, shown in Fig. \ref{U_excited_DeltaT}. As we increase the delay the synchronization
level decreases significatively, the global excitement of the network is
lower and the feedback input is lower. Because of the frequency
dependence of the response of the FHN equation on the input, the frequency of
oscillation decreases (up to $20\%$) with increasing delay. Following this model we would
expect that two equal networks, strongly connected through paths with long conduction delays, would oscillate with lower frequencies and
less synchrony than those connected with short conduction delays.

This is indeed found in the brain, where the oscillatory activity of
connected networks oscillates with lower frequencies in the case of
connections that go along longer
distances. For instance while local sensory integration evolves with a fast gamma (25-70 Hz) dynamics,
multisensory integration evolves with an intermediate beta (12-18Hz) dynamics,
and long-range integration during top-down processing evolves with a temporal
dynamics in a low theta/alpha (4-12Hz) frequency range \cite{Stein}. Many authors have studied the role of time-delays and nonlocal excitatory
coupling in the generation of synchronous rhythm in the brain, although the detailed
mechanisms behind this fact are still under investigation \cite{Kopell,Er,Rub}.

\clearpage
\newpage
\subsection{Periodic force and feedback, reinforcement of an external stimulus.}\label{2-4}

If we study the system excited by noise with $D=0.005$,the
feedback term described above with $A=0.5$ and an external periodic force $I(t)=A\cos(2\pi f t)$ with
frequency $f=0.55$ (period $T=1/f$), we observe that after the
transient time the density moves along the attractor in the phase space plane in a periodic
way, with the period of the exciting force. We can see that two different
behaviors coexist in a single network:  an almost Gaussian distribution 
is moving along the attractor trajectory, whereas a secondary maximum is
stable and localized at the non-firing position.  The fraction of neurons above the
threshold, $n=\int_0^{\infty}\rho (u,v,t-\Delta T) du
 dv$, oscillates now between $0\%$ and $90\%$, the remaining $10\%$ is 
 at rest whereas without the feedback term the maximum $n$
 is $60\%$. In this case the feedback term and the noise helps the network to
 synchronize giving a very strong response to a signal, in other words:
 such feedback structures helps the system in amplifying a signal. 

These recurrent neuronal circuits where synaptic output is fed back
as part of the input stream, exist naturally in cerebral structures as
amplifier structures. Neural connectivity is highly recurrent (for instance
Layer IV is believed to be the cortical signal amplifier of thalamic signals).


\begin{figure}[h]
\begin{center}
\epsfig{figure=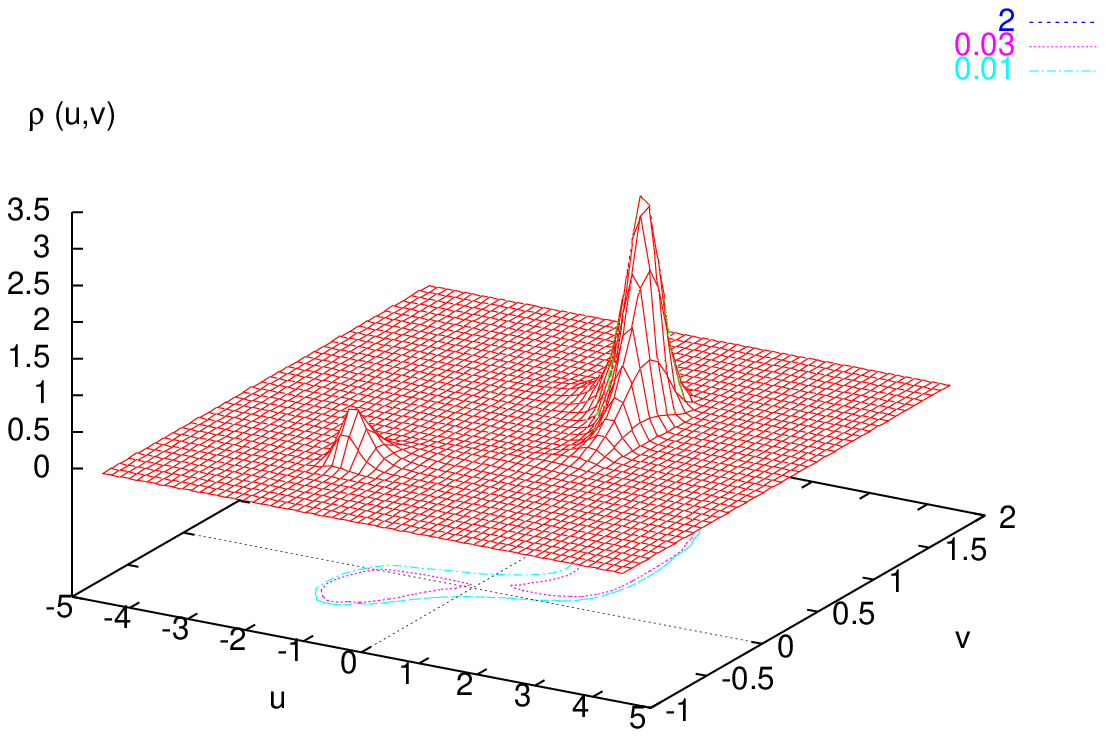,height=7cm,width=6.5cm}
\epsfig{figure=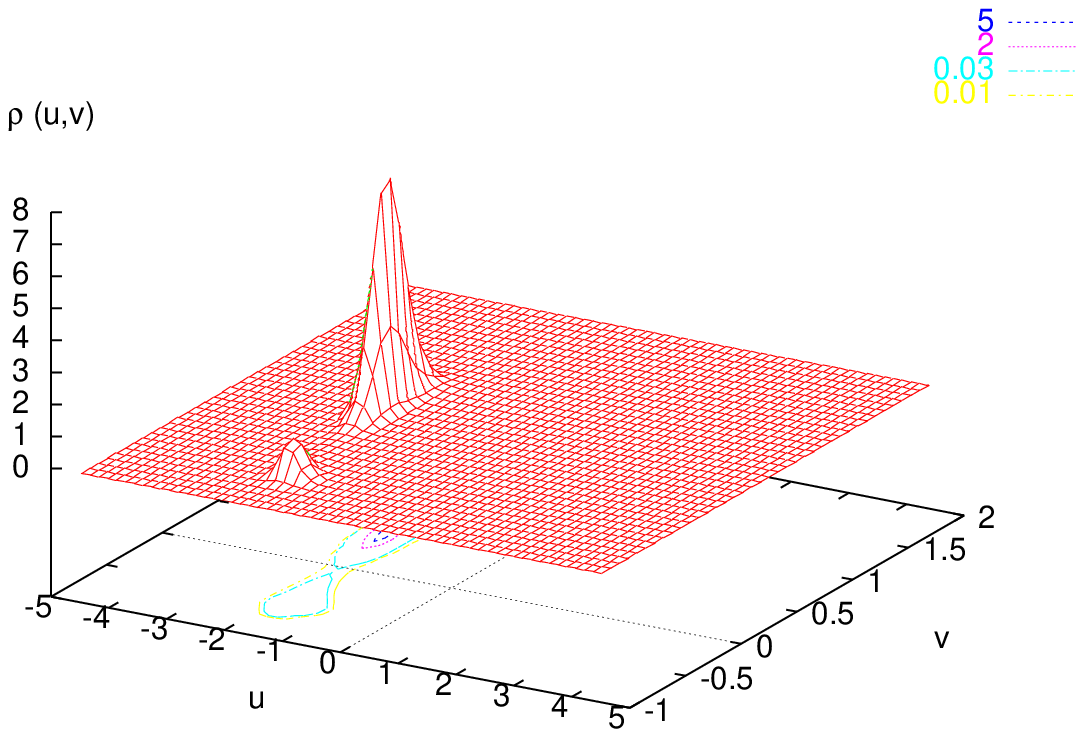,height=7cm,width=6.5cm}
\caption{\label{rho_excited_DeltaT} Density for feedback with $\Delta T=3/4T=3.00$ and
  periodic excitation with $A=0.5$, shown at
 at times $t=35.84$ and $t=40,96$. The network is strongly synchronized, all
 the neurons evolve in a similar state.}
\end{center}
\end{figure}

\section{Conclusion}

The Fokker-Planck equation has been used to explore the time evolution of the probability
density of a neural network excited by noise. Integrating numerically this equation
one can obtain a quick overview of the system dynamics, its moments and
spectral density, the attractors and most probable states of the neuron
ensemble. We can explore different basic mechanisms that show new emergent properties depending on
the connectivity of the network and the noise intensity. 

An intermediate level of noise alone, combined with nonlocal excitatory
interactions, can give rise to coherent, strongly synchronized, oscillations.
There are several competing mechanisms: noise, dispersion, nonlinearity and
nonlocal interactions (the feedback term) \cite{LVKonotop}. The balanced combination of all of
them leads to the formation of a global oscillatory and robust state. 

This dynamics also helps the system to follow any external weak excitation, such as a periodic force. If the
 parameters are properly chosen as much as $95\%$ of the neural population can
 be synchronized improving the sensibility of the system to weak external signals. If the delay time or if the noise is out of
 the adequate range the global oscillations process does not take place.

\subsection{Acknowledgement.}
L. V\'{a}zquez wants to thank the partial support of the European Project COSYC of
SENS (hprn-ct-00158).

\end{document}